# $M|G|\infty$ Queue Busy Cycle Renewal Function

Manuel Alberto M. Ferreira

manuel.ferreira@iscte.pt

Instituto Universitário de Lisboa (ISCTE-IUL), ISTAR-IUL, Lisboa,

PORTUGAL

**Abstract**

We present formulas to compute the busy cycle renewal function for the $M|G|\infty$ queue and exemplify for some service time distributions. The busy cycle renewal function value in *t* is the number of busy periods that begin in $[0, t]$. This number is of crucial importance in practical applications of $M|G|\infty$ queue.

**1. Introduction**

In $M|G|\infty$ queue, $\lambda$ is the customers' Poisson process arrivals rate, in the moment it arrives at the system begins each customer' service, which has a positive distribution, which distribution function we call $G(\cdot)$, and mean $\alpha$. So, $\alpha = \int_0^\infty [1 - G(v)]dv$. The traffic intensity is given by $\rho = \lambda\alpha$.

In a queuing system, we call busy period a period that begins when a customer arrives there and find it empty; it ends when a customer leaves the system letting it empty; and there is always at least one customer present. So, in a queuing system, there is always a sequence of idle periods and busy periods.

Let's consider the $M|G|\infty$ system with time origin at the beginning of a busy period. The instants $0, t_1, t_2, ...$, at which a busy period begins, are the arrival epochs of a renewal process (Takács,1962). We say that a cycle is complete when a renewal occurs, that is: a busy period begins. These cycles are busy cycles, ( ) and their time is a positive random variable that we will call *Z*.

So,

$$Z = B + I \quad (1.1),$$

Where *B* is the busy period time and *I* the idle period one.

In (Takács,1962) the author proved that *B* and *I* are stochastically independent and the Z Laplace transform, $\bar{Z}(s)$, is

$$\bar{Z}(s) = 1 - \frac{1}{(s + \lambda)P_{00}(s)} \quad (1.2),$$

Where $P_{00}(s)$ is the Laplace transform of $p_{00}(t) = e^{-\lambda \int_0^t [1-G(v)]dv}$, the probability of the system emptiness at *t*, being initially empty.

Takács also showed that

$$E[Z] = \frac{e^\rho}{\lambda} \quad (1.3)$$

And

$$E[Z^2] = 2\lambda^{-1}e^{2\rho} \int_0^\infty \left[e^{-\lambda \int_0^t [1-G(v)]dv} - e^{-\rho}\right] dt + 2\lambda^{-2}e^\rho \quad (1.4).$$

Being $I$ exponentially distributed with parameter $\lambda$, as it happens with any queue system with Poisson arrivals, its Laplace transform is $\bar{I}(s) = \frac{\lambda}{\lambda+s}$ and the ratio $\frac{\bar{Z}(s)}{\bar{I}(s)}$ gives the expression $\bar{B}(s) = 1 + \frac{1}{\lambda}\left[s - \frac{1}{P_{00}(s)}\right]$ for the $B$ Laplace Transform, (Stadje, 1985) which inversion is a complex task, except for some service time distributions, see (Ferreira, 1991, 1995, and 1998).

Our work will focus on the $M|G|\infty$ queue busy cycle, in its renewal function study.

## 2. The $M|G|\infty$ queue busy cycle renewal function

The renewal function, $R$, of a renewal process is $R = 1 + F + F^{*2} + F^{*3} + \cdots$, where $F^{*n}$ is the $n$-th interrenewal time convolution with itself distribution function, see, for instance (Çinlar, 1975). It gives the mean number of renewals in $[0, t)]$.

For instance, in the application of this model to unemployment situations, a busy period is a period of unemployment. And, in illness situations, a busy period is an epidemic or pandemic period. See, about this kind of applications (Ferreira, 2003 and 2003a).

To compute the $M|G|\infty$ queue busy cycle renewal function, using the Laplace transform:

$$\bar{R}(s) = \frac{1}{s} + \frac{1}{s}\bar{Z}(s) + \frac{1}{s}\bar{Z}(s)^2 + \cdots + \frac{1}{s}\bar{Z}(s)^n + \cdots = \frac{s^{-1}}{1-\bar{Z}(s)} = \frac{s^{-1}}{1-1+\frac{1}{(s+\lambda)P_{00}(s)}} =$$

$$\frac{(s+\lambda)P_{00}(s)}{s} = P_{00}(s) + \lambda \frac{1}{s} P_{00}(s). \text{ So, } R(t) = p_{00}(t) + \lambda \int_0^t p_{00}(u) du \text{ and}$$

$$R(t) = e^{-\lambda \int_0^t [1-G(v)]dv} + \lambda \int_0^t e^{-\lambda \int_0^u [1-G(v)]dv} du \quad (2.1)$$

Note that

$$\lim_{t\to\infty}\left[R(t) - \frac{\lambda}{e^\rho}t\right] = \lim_{t\to\infty}\left[e^{-\lambda \int_0^t [1-G(v)]dv} + \lambda e^{-\rho}\int_0^t e^{\rho-\lambda \int_0^u [1-G(v)]dv} du - \frac{\lambda}{e^\rho}t\right]$$

$$= e^{-\rho}$$

$$+ \lambda e^{-\rho} \lim_{t\to\infty} \int_0^t \left[e^{\rho-\lambda \int_0^u [1-G(v)]dv} - 1\right] du = e^{-\rho}$$

$$+ \lambda e^{-\rho} \int_0^\infty \left[e^{\lambda \int_u^\infty [1-G(v)]dv} - 1\right] du.$$

So, we can conclude that

$$\lim_{t\to\infty}\left[R(t) - \frac{t}{E[Z]}\right] = \frac{VAR[Z] + E^2[Z]}{2E^2[Z]} \quad (2.2),$$

As it as to be with any renewal function.

-As $e^{-\rho} \leq p_{00}(t) \leq 1$,

$$p_{00}(t) + \lambda e^{-\rho} t \leq R(t) \leq p_{00}(t) + \lambda t \quad (2.3),$$

And still

$$e^{-\rho}(1 + \lambda t) \leq R(t) \leq 1 + \lambda t \quad (2.4),$$

We conclude that

$$\lim_{\alpha \to 0} R(t) = 1 + \lambda t \quad (2.5).$$

As it must be because, when the service time is null, when it arrives each customer begins a busy period. And the arrival instants, in the $M|G|\infty$ system occur according to a Poisson process at rate $\lambda$.

-

$$\frac{d}{dt} R(t) = -\lambda p_{00}(t)(1 - G(t)) + \lambda p_{00}(t) = \lambda G(t) p_{00}(t) \geq 0$$

So $R(t)$ increases with $t$ (evidently).

**3. $R(t)$ expressions for some service time distributions**

Now we will present $R(t)$ expressions for various service time distributions, using formula (2.1). For some cases only bounds were possible to determine:

1. $G(t) = 1 - \frac{(1-e^{-\rho})(\lambda+\beta)}{\lambda e^{-\rho}(e^{-(\lambda+\beta)t}-1)+\lambda}, t \geq 0, -\lambda \leq \beta \leq \frac{\lambda}{e^\rho - 1}$, see (Ferreira,1998)

$$R(t) = e^{-\rho}(1 + \lambda t) + (1 - e^{-\rho})\frac{\beta}{\lambda + \beta} e^{-(\lambda+\beta)t} + (1 - e^{-\rho})\frac{\lambda}{\lambda + \beta}, -\lambda < \beta < \frac{\lambda}{e^\rho - 1} \quad (3.1)$$

2. **Deterministic service with value $\alpha$:** $G(t) = \begin{cases} 0, & t < \alpha \\ 1, & t \geq \alpha \end{cases}$

$$R^D(t) = \begin{cases} 1, & t < \alpha \\ 1 + \lambda e^{-\rho}(t - \alpha), & t \geq \alpha \end{cases} \quad (3.2)$$

3. **Power function service with parameter c, $c > 1$ (for $c = 1$ this is the uniform distribution in $[0, 1]$):** $G(t) = \begin{cases} t^c, & 0 \leq t < 1 \\ 1, & t \geq 1 \end{cases}$

$$R^P(t) = \begin{cases} e^{-\lambda\left(t-\frac{t^{c+1}}{c+1}\right)} + \lambda \int_0^t e^{-\lambda\left(u-\frac{u^{c+1}}{c+1}\right)} du, t \leq 1 \\ e^{-\rho} + \lambda \int_0^1 e^{-\lambda\left(u-\frac{u^{c+1}}{c+1}\right)} du + \lambda e^{-\rho}(t-1), \ t > 1 \end{cases} \quad (3.3)$$

4. **Exponential service times:** $G(t) = 1 - e^{-\frac{t}{\alpha}}, t \geq 0$

$$R^M(t) = e^{-\rho\left(1-e^{-\frac{t}{\alpha}}\right)} + \lambda e^{-\rho} \int_0^t e^{\rho e^{-\frac{u}{\alpha}}} du \quad (3.4)$$

5. **Service NBUE-New Better than Used in Expectation**

If the service distribution is NBUE, with mean $\alpha, \int_t^\infty (1 - G(v))dv \leq \int_t^\infty e^{-\frac{v}{\alpha}} dv$, see (Ross,1983). Then

$$R^{NBUE}(t) \leq R^M(t) \quad (3.5)$$

6. **Service NWUE-New Worse than Used in Expectation**

If the service distribution is NWUE, with mean $\alpha, \int_t^\infty (1 - G(v))dv \geq \int_t^\infty e^{-\frac{v}{\alpha}} dv$, see (Ross,1983). Then

$$R^{NWUE}(t) \geq R^M(t) \quad (3.6)$$

7. **Serviço DFR-Decreasing Failure Rate**

If the service distribution is DFR, $1 - G(t) \geq e^{-\frac{t}{\alpha}\frac{\gamma_s^2}{2}+\frac{1}{2}}$, being $\gamma_s$ the service coefficient of variation, see (Ross,1983). Then

$$R^{DFR}(t) \leq e^{-\rho\left(1-e^{-\frac{t}{\alpha}}\right)e^{\frac{1-\gamma_s^2}{2}}} + \lambda e^{-\rho e^{\frac{1-\gamma_s^2}{2}}} \int_0^t e^{\rho e^{\frac{1-\gamma_s^2}{2}-\frac{u}{\alpha}}} du \quad (3.7)$$

8. **Serviço IMRL-Increasing Mean Residual Life**

If the service distribution is IMRL $1 - \frac{\int_0^t [1-G(v)]dv}{\alpha} \geq e^{-\frac{2t\alpha}{\mu_2}-\frac{2\alpha}{3\mu_2^2}\mu_3+1}$ being $\mu_r$ the $r$ order moments around the origin, see (Brown, 1981) and (Cox, 1962). Then

$$R^{IMRL}(t) \geq e^{-\rho}\left(1 - e^{-\frac{2t\alpha}{\mu_2}-\frac{2\alpha}{3\mu_2^2}\mu_3+1}\right) + \lambda e^{-\rho} \int_0^t e^{\rho e^{-\frac{2u\alpha}{\mu_2}-\frac{2\alpha}{3\mu_2^2}\mu_3+1}} du \quad (3.8)$$

**Note**:

-The last four service distributions presented, and related to the exponential, are widely used in reliability theory.

**4. Conclusions**

After a short study about $R(t)$ and its properties, we present, in the situations of some service time distributions, formulas to its evaluation. Deserve to be highlighted the cases of the distributions NBUE, NWUE, DFR and IMRL widely used in reliability theory, and important in practical applications.